\documentclass[11pt,reqno]{amsart}
\usepackage[top=1.25in,bottom=1.25in,left=1in,right=1in]{geometry}
\usepackage{amsfonts, amsmath, amssymb, amscd, amsthm, graphicx}
\usepackage{enumitem}
\usepackage[final, 
colorlinks=true]{hyperref}
\usepackage{verbatim}
\usepackage{url}

\usepackage{tikz,tikz-cd,amscd}
\usetikzlibrary{calc}
\usepackage{tikz,mathdots,color}
\usetikzlibrary{decorations.pathreplacing,decorations.markings}
\usetikzlibrary{shapes,shapes.arrows,positioning}

\tikzset{->-/.style={decoration={
  markings,
  mark=at position #1 with {\arrow{>}}},postaction={decorate}}}
\tikzset{->-/.default=.5}

\theoremstyle{plain}

\newtheorem*{claim*}{Claim}

\theoremstyle{definition}

\numberwithin{equation}{section}
\numberwithin{theorem}{section}

\newcommand{\param}%
	{{\mathchoice{\mkern1mu\mbox{\raise2.2pt\hbox{$\centerdot$}}\mkern1mu}%
	{\mkern1mu\mbox{\raise2.2pt\hbox{$\centerdot$}}\mkern1mu}%
	{\mkern1.5mu\centerdot\mkern1.5mu}{\mkern1.5mu\centerdot\mkern1.5mu}}}

\begin{document}

\title{Navigating the AI crisis: a humble guide for students}


\author{Tarik Aougab}
\address{Department of Mathematics, Haverford College, Haverford, PA 19041, USA}
\curraddr{}
\email{taougab@haverford.edu}

\begin{abstract} We present some basic recommendations for (mostly graduate) students when it comes to choosing to use, or avoid, AI in their mathematical research and education. These include: maintaining a connection to the tried and true ways for learning mathematics; focusing on communication and the sharing of ideas with other people; thinking critically about the politics of labor and organizing with your cohort and others at your home institution; and being in touch with the realities of the impact of AI infrastructure on land and people. 
\end{abstract}

\date{\today}

\maketitle 

The impact of AI on professional mathematics was a central focus at the most recent International Congress of Mathematics (ICM) in Philadelphia. Keynote addresses and major panel discussions were dedicated to thinking through the inevitable changes that are coming to academic mathematics due to the advent of models that can assist in solving major conjectures. And in the last few weeks, \href{https://proofsandprompts.com/}{discussion regarding AI in mathematics has exploded} (see also \cite{ar:Weinreich}, \cite{ar:Hofer}, \cite{ar:Boldt}). While I definitely disagree with \href{https://sites.google.com/view/icm2026boycott/home}{many of the choices made by the IMU regarding this ICM}, this wasn’t one of them; mathematicians \textit{must} think deeply about the role that AI assistance will play in our discipline. 

The only way to do this is to reflect on \textit{our own} roles and goals. For example, \href{https://www.youtube.com/watch?v=M0--ZH1lOzg}{Terry Tao’s ICM address} \cite{Talk:Tao}– geared towards constructing a set of recommendations for the creation and digestion of pure mathematics in the time of AI– was organized around first being clear about what the hell it is that mathematicians are up to and what we want to achieve. 

\section{What the hell are we up to and what do we want to achieve as mathematicians?}

There are many possible answers, for instance: generating absolute truths, increasing popular interest in mathematics, deepening our collective understanding of mathematical structures, solving major open problems, posing interesting conjectures, cultivating human beings who enjoy problem solving and theory building, developing mathematical theory that leads to the invention of new technologies that benefit life, creating an environment for which the study of mathematics is encouraged (or even joyful), creating quantitative tools that serve everyday people and promote justice, grow our understanding and in turn our appreciation of nature, land, and ourselves, etc.

Up until now, research mathematics was structured around the premise that success in one of these arenas, namely resolving the final step in a major open problem, should dictate which mathematicians get to evaluate the importance of all of these possible targets. This was justified by the assumption that such targets were pairwise positively correlated. Putting aside whether or not this was ever true, it seems likely that in the near to distant future, it no longer will be. For example, AI is capable of asserting new truths while contributing near absolute zero to the deepening of our collective understanding. It is also capable of increasing popular interest in mathematics while simultaneously jeopardizing some of the most conventional and well trodden pathways for studying and doing math. The number of proclamations asserting the ``end of human mathematics'' (it’s worth pointing out that CEOs and other employees of AI start-ups are vastly over-represented in the group of people making such claims) grows by the day.

The premise for this piece is simple: no matter what one thinks the central goals or purposes of mathematics are, the achievement of those goals by humans is impossible without ensuring that human beings have ample access and opportunity to study and think about mathematics. Of course, there are certainly those in our community who believe that the goal that matters more than any other is the pursuit of truth. So as long as someone, or some\textit{thing}, is generating mathematical truth, mathematics can flourish. And in fact – so this line of reasoning goes– it is the responsibility of the mathematical community to empower those people and/or inanimate objects that are “best”\footnote{In quotes to stress that this is obviously not well-defined.} at pursuing those truths. If X equals the total amount of resources that a society puts aside for the production and practice of mathematics, a reasonable implication of this thought process is that X should be divided up among the set of all nouns in accordance with the capacity of each such noun– be it a computational model, a human, a plank of wood, whatever–  to produce mathematical truth. 

To me, and thankfully to many of the other mathematicians who’ve weighed in on this conversation before (see \href{https://proofsandprompts.com/}{the many wonderful posts here} \cite{Blog:PAP}), mathematics is by definition a human practice. It is therefore dead the moment there does not exist human beings doing it. To fully embrace this point of view, it’s important to flesh out what it really means to say that ``mathematics is a human practice'' because that phrase has become a bit of a slogan, and I want to be crystal clear. My perspective, which I think is pretty uncontroversial, is that mathematics evolved organically and over time alongside many other human activities. As in, it is not as if one day a collection of people got together and consciously decided that mathematics was worth doing and ever since then we’ve done mathematics, in the same way that there was almost surely never a singular point in human history where people got together and decided that we should dance, or play music, or make art, or cultivate distinct cuisines and languages, or develop architecture, or hone analytical methods of thought.

All of these things came into being over time because they are deeply intertwined with the experience of being human. They developed because doing them brought joy to our lives and deepened our connection and understanding of our surroundings and of one another. And I submit that no matter how emphatically you might claim otherwise in the broader impacts section of your next grant, \textit{you} don’t really do mathematics because it is useful for the world or because it will or may lead to some such application 100 or 1000 years down the road (even though it absolutely might). You do it because you love doing it, and because doing it allows you to tap into the millennia-long tradition of other humans doing it because \textit{they} loved it. Mathematics facilitates the human appreciation of the natural world, including of our own minds. 

So I hope you’re willing to accept the claim that mathematics requires humans, and conversely, \textit{humans require mathematics}, not simply for the ways in which it eventually betters our lives through technological advancement but actually for its practice as a thing unto itself. A future in which all of mathematics is being done by machines is therefore not really possible, in the same way that a future that contains humans and in which all speaking or all dancing is being done by machines is impossible. On the other hand, a future in which it is virtually impossible for a human to pursue mathematical expertise at a professional level \textit{could} come to pass.

\section{What the hell is the author up to and what does he hope to achieve in writing this? 
}

It is easy for a distinguished mathematician near the end of their career to entertain the idea of there being no such thing as professional mathematics. I’ve seen several senior people conceiving of such a future as being overall net good for society so long as it coincides with AI churning out more math than humans ever could. My main inspiration for writing this piece is feeling frustrated at how little consideration is given to young researchers and students in such fantasies. Our students wish– just as every generation of scholars wished before them– to find their niche and contribute meaningfully to the collective body of knowledge. As AI optimists frequently point out, there is nothing inherent or good about structuring our community’s evaluative process around the generation of proofs or the solving of conjectures. And so, as the argument goes, perhaps this new technology will free us to start valuing the other, arguably far more important contributions to mathematical knowledge production besides being the last person in a very long line of other people to add an insight that finally cracks open an old problem. 

This is all well and good, but what \textit{is} inherent to the human practice of mathematics is the opportunity for each person to make the material their own, or to reshape and reframe it in ways that make sense for their own minds and for the minds of other people who may think like they do. People need to be given the time, space, and resources to do this and in general to think deeply about mathematics. \textit{This} is what we need to preserve; the concept of the research mathematician-as-profession is just a proxy for this. For example, perhaps it wouldn’t be so bad for the career to die if its death coincided with comprehensive economic policies– like a genuine universal basic income– whose existence permits people to spend their time thinking, writing, and creating their own mathematics. A future in which this time is afforded to even fewer humans than it is today is a terribly bleak one. If we hope to avoid that future, we need to start thinking not only about how (or if) we integrate AI into our research workflows, but also about how to guide our youngest and newest community members through this turbulent time. Without them, there is actually by definition no future for mathematics.

So this is what I’d like to focus on here. In particular, this piece is dedicated to offering some recommendations and suggestions for current and future students in the mathematical sciences. How do students navigate this moment in our discipline’s history? What should or shouldn’t they be doing to increase the likelihood of success in the professional mathematics community? How are students supposed to learn when such a significant proportion of their role models are spending more time feeling sorry for themselves than they are mentoring the next generation? 

Too often, if students are spoken to at all in a format like this, it only happens as an afterthought. I’d like to flip that script: in what follows, I’ll address students directly with a special focus on graduate students. To the extent that there are messages directed towards senior researchers, they will be implicit and it might be necessary to read between the lines to find them. Some of the suggestions are straightforward and uncontroversial; others are broader and more expressly political in nature. I feel that (some version of) all of them would have made sense as advice to students even without the existence of AI, but now they are more urgent than ever. They are not meant to be exhaustive and you, the reader, have no obligation to agree with all of them. 

Finally, prescription is intentionally avoided; the goal is neither to convince you to use AI or to avoid it completely. Besides being a practicing mathematician and caring about students, I have no particular expertise.  I'll have succeeded if by writing and circulating this, students feel less ignored in this discussion and more empowered to take ownership over the various decisions that await them.

\section{Divide your time between studying the fundamentals and what may become the \textit{new} fundamentals. }

I can't remember a time when I actually needed to compute the solution of an ODE exactly as part of my research. I can admit that practice computing such solutions was a helpful part of my mathematical development. By the turn of the last century, technology and mathematical culture had evolved to such a point that many calculations which at one point might have constituted on its own a reasonable contribution to the field had been made obsolete by computer algebra systems. No one ever seriously took this as a sign to stop teaching the art of e.g. solving integrals, and nor should they have. At the same time, offloading some of these calculations onto machines freed up space in our psyches to start asking different kinds of questions.  

This development is not unique to physical technological advancement; the history of mathematics is replete with times when \textit{ideological} advancements turned what used to be very challenging problems into simpler ones and in the process, opened the door to new flavors of curiosity. When this happens, there is often a process of recalibration in the field, during which the size, complexity, and perhaps abstraction of a fundamental building block of argumentation increases, or at least shifts in some way. 

The use of AI in mathematics represents both types of advancement. It is obviously a new physical technology, and that technology may perhaps be used to generate ideological breakthroughs. I’ve heard from some students that their advisors are encouraging them to stop studying and reading and practicing mathematics the old way and to adapt completely to whatever becomes the new building blocks of mathematical reasoning. This is just as misguided as it would have been to stop teaching integration when Wolfram released Mathematica. Your adviser knows a lot of mathematics, but they do not know anymore than you about how to learn mathematics in a time when AI exists. 

The advice here is actually kind of simple: students of mathematics should do what they can to keep up with the trajectory of their chosen field, while simultaneously following best practices for how to learn. In the past, this meant reading the standard texts, struggling with exercises both on your own and in groups with your cohort, and subscribing to the appropriate ArXiv tags. It still does mean all of that! But perhaps it also means being in conversation with your peers and colleagues about the role AI is playing in your subfield, regardless of whether you decide to use AI or not. What types of problems has it solved? What does it seem to be good at versus what can it still not accomplish easily? This involves asking yourself questions not only about AI itself, but about what others in your field think about AI. As it stands, this is a community of humans, and humans are still the ones who make decisions about resource allocation, hirings, and performance evaluation (and in a future where such decisions are no longer human-made, it will still be humans who choose to automate those decisions and how to automate them). So for example, you might wonder which sorts of computations it would be reasonable for an expert to expect AI to be able to execute on its own. Since it’s plausible that such work will be heavily devalued in the near term, having a solid feel for this will help you to calibrate what your colleagues will and won’t find mathematically interesting when it comes to judging your work. 

\section{Focus on developing your skills of exposition and communication.
}

As stressed above, doing mathematics is fundamentally about sharing our understanding of structure and pattern with one another. It is therefore a  failure of professional mathematical culture that as a community, we have never \textit{really} strayed from the opinion of Hardy that ``exposition is work for second-rate minds.'' Of course in reality, when an important new theorem is proved, its true impact can only be realized after a tremendous amount of hard thinking, re-interpretation, and writing and speaking about it by many many people. 

And yet, it is sadly still par for the course for top tier graduate programs to spend little or no effort training students in the mathematical art of communication. This has set up a tragic dynamic which has allowed for large language models to claim an outsized influence over the trajectory of our discipline. As an example of what I mean, consider MathOverflow. The idea behind it is wonderful: a place where people can virtually gather to discuss research-level mathematics questions and where genuine connections can be made between question-asker and perhaps many different question-answerers.

Unfortunately, its reputation among many – perhaps most of all the younger and less established community members– is that it is an unbearably intimidating hellscape where you’re just as if not more likely to be chastised for asking something too basic as you are to be graced with an actual answer. Who in their right mind would choose to subject themselves to that sort of treatment when they could just ask ChatGPT? By de-emphasizing the essential social skills that mathematics depends upon, we’ve made it far less likely that any one of us would pick the messiness of human-to-human interaction over the easy and sycophantic– but also not insulting or scary– LLM. 

As part of the next generation of mathematical scholars, you can help the community to shift course. This would be admirable and extremely important even if AI didn’t exist. But because it does, it is absolutely necessary for the survival of our community. There may come a time in the near future when most research-level proof generation is being done by AI models. In that sort of landscape, the work of human experts becomes what it already was for a fraction of $1-\epsilon$ of every single human mathematician’s working hours (even those among us who’ve proved landmark theorems spend most of their time not proving such theorems): interpretation, communication, re-framing, remixing, etc. 

``\textit{But what if AI can do all of that better than humans too?}'', I hear you ask. By definition, AI cannot be better than humans in the activity of human beings talking to one another. So insofar as humans sharing ideas with other humans is an essential part of what it means for mathematics to exist, this question is actually beside the point. In fact, there is nothing automatically wrong or bad about a human using an AI to help themselves understand some piece of mathematics in order to then go and explain it in a new or interesting way to another human or group of humans. But this is only doable at all if human mathematicians learn to value– and even reward – conversation.

To prepare yourselves for these changes, I recommend practice: organize student seminars where people can speak about their own work in a judgment-free environment, talk to one another informally and as much as you can about what you’re learning, and write up original lecture notes from a class and make them available to others. Seek out some of the best expository writing in your subfield and see if you can pick out what makes it as good as it is. And, make yourself known in your mathematical niche as a person who loves talking to and sharing with others.

\section{Talk to your fellow students about working conditions across all levels of the program.  
}

One of the irrefutable lessons that came out of the crackdown around the anti-Zionist, anti-genocide student movement is that the transfer of power at the university from academics to board members and governmental functionaries, is essentially complete. On campus, the desire of a handful of donors to silence a particular conversation evidently outweighs the pursuit of truth, ethics, or really any other virtue that the university might claim to care about in its fancy Latin motto. Schools don’t have to be organized in that way– and in fact they haven’t always been– but we can’t act as though they aren’t in fact like this currently.

What happens when decisions on campus are made by business people instead of students, researchers, and teachers? The overall quality of instruction and the long term, hard-to-quantify benefits of a successful university education are sidelined. Increasingly, your academic employer will try to squeeze as much out of you as it possibly can, while providing as little as it possibly can to you in the form of mentorship, training, and of course money. This has absolutely nothing to do with the personalities or politics of your direct superiors, such as your adviser or the department chair. We are talking instead about the dynamics of an entire labor market.

Recently, the Trump administration \href{https://www.whitehouse.gov/releases/2026/07/45502/}{announced its intention} to siphon even more money away from traditional universities and towards privatized AI development. So these dynamics will most likely only be furthered by this current flurry of activity around AI. It is important to keep in mind: the capital that has been mobilized in order to bet on the potential of AI to generate unprecedentedly large returns represents the largest in modern human history. So even in a world where AI does not end up living up to the hype, its impact will be seismic on the university if for no other reason than the fact that in such a future, recessions that would make the housing crash of 2008 look like a blip by comparison are likely. 

Older professors tend to like to wax romantically about graduate school: very few responsibilities and endless time to just think about hard mathematics. Your generation has no such luxury because you (very literally) will not be able to afford that level of oblivion. Talk to one another about your working conditions; make it public knowledge in the department what everyone is receiving in compensation and what is expected of each person in the program. Be aware of how decision-making at all levels of the university (the way your adviser allocates their grant money; the way the department budgets the money allotted to it by the dean of sciences; the way the provost distributes resources to the various academic departments; the way the university chooses to spend its money on academic development and existing infrastructure– like fixing up a dilapidated dormitory so that students actually have a safe place to live– as opposed to building a new stadium; the international partnerships your university pursues or stays away from; the ways your school treats the land on which it is based and the relationship it maintains with surrounding neighborhoods and communities; and the way that the state and federal government funnels money to or away from the institution) impact the material conditions of your job. In fact, I encourage students to dedicate a seminar– in addition to the one recommended in the previous section – to discussing these issues. 

One way to think about this is that you have it harder than your predecessors. And that’s in many ways true. But another way to think about it is that your predecessors are (on average) politically stupid because they never had to think seriously about these matters. The economics of AI is all about labor; indeed, the capacity of AI to lower labor costs (economists’ preferred euphemism for firing people and replacing them with machines) is what accounts for its obscene valuation. For these reasons, even the most apolitical among us will need to smarten up. The next crop of mathematicians must be populated by people who understand the importance of labor power.

\section{That last point generalizes even further to students in other departments and fields. 
}

The capacity of the modern university to exploit you is predicated on its ability to isolate you as much as possible from your peers. Administrators argue frequently that it’s just common sense to treat workers in the history department completely differently to those in the chemistry department, and they point to the obvious differences in the two jobs to back up this claim. A chemist’s daily life is on the surface totally different from that of a historian’s, so do the administrators have a point?

It’s actually sort of amazing how selective the logic of an administrator or board member can be. On the one hand, they can go on for hours about new tools of financialization or clever abstractions cooked up by the legal team for avoiding taxes on land holdings; anything is possible with the right mindset and enough creativity! But the idea of an institution-wide employee contract is simply beyond the pale because some people synthesize chemicals and others write books. Of course, there is such a thing as an industrial union. At an airline, it would protect the basic employment rights of people with very different jobs, like ticket agents, baggage handlers, and mechanics.

What does any of this have to do with AI, or with mathematics specifically? Long before the modern LLM was even a sparkle in anyone’s eye, efficiency-obsessed administrators have been trying to gut language, history, philosophy, music, and cultural and ethnic studies departments on the basis that they cost more than they’re ``worth.'' Here, ``worth'' is in quotes because the administrator doesn’t know how to measure the worth of developing a generally informed and well-rounded population, and so the ``worth'' of the history department instead ends up being quantified in terms of, for instance, the number of dollars that history major alumni are in a position to donate to the university. Of course, the more resources the administrator takes away from the department and the more the university prioritizes the departments that score well on this or some other metric whose use is justified only by the ease with which it is computable, the fewer students choose to major in those neglected departments, the more the administrator feels justified in taking even more away, and so on and so forth. 

For decades, mathematics departments have felt invulnerable to this vicious cycle. We could claim that what we teach and research has tremendous utility, and students majoring with us could likely go on to achieve incredible levels of financial success. Now that it is possible that an AI can do almost anything we teach our students to do, this claim no longer holds the water it once did. The problem here is that we should have never bought into the stupid metrics that the administrator used to rank us so much higher than our humanities colleagues in the first place, and we ignored the cries of those colleagues over the last four decades to our peril. 

They had it right all along: the purpose of teaching history is not to maximize the success levels of students in the labor market, nor is it to maximize the number of people who make enough to get to exit that market and can then begin to buy and sell other people’s time. It’s to do something else entirely, not least of which would be to help develop a society that can sense when it is being robbed blind by Silicon Valley. Likewise, mathematicians need to internalize what the real utility is for teaching and thinking about mathematics. 

One of the best ways to do this is to form connections with your colleagues in other departments who most likely have years more practice thinking about this than you do. And after doing so, you’ll hopefully be in a much better position to organize for whatever it is that you all want and need on campus. If students at your institution are already organized into a union, that is a great place to start. If not, help lay the groundwork for one.

You might be thinking: \textit{with AI advancing, machines can very easily take our place and so we have less leverage than any generation of students before us ever did.} This is a fair point, but there have been many times throughout history during which the greatest gains for workers occur in reaction to the most dramatic accumulations of wealth and inequality. We cannot predict all – or even really any– of the ways that the advent of AI will impact labor in the long run, but whatever is to come, it can’t hurt to face it together. 

\section{Take the time to learn about the place where you are and the people who live there. 
}

If you’re plugged into the conversations in the media around AI, you’re probably finding it difficult to avoid the recent onslaught of hot takes from AI industry people claiming that the massive public opposition to data centers is ill-informed, uneducated, and probably also somehow astroturfed. These people– many of whom have glaring conflicts of interest as they are either massive investors in AI infrastructure or are employed by such capitalists– argue that people are just being emotional and that this wave of anti-AI sentiment is the same sort of uneducated fear of change that accompanies any major technological advancement.

These are also the people who, as a group, do an absolutely abysmal job of what's recommended in this section. They do not understand what people actually need or want, in part because they have spent years convincing one another that they have the right and the power to \textit{dictate} what people need and want. Evidently, if someone living near a data center resents it because they no longer have access to fresh water or they can’t get enough sleep because of the tens of decibels it constantly outputs, that person just doesn’t know yet how much they ``want'' the center. 

\href{https://arxiv.org/abs/2608.02859}{In a recent essay}, Max Weinreich argues (I believe correctly) that this anti-AI surge is genuine and that it is based in real-world, serious concerns and problems that people have with the technology. And if mathematicians don’t think extremely intentionally about how we engage with it, we run the risk of guilt by association. It would be hard to imagine a worse way of injuring our public credibility. Even though there is much in Tao’s ICM address I can get behind, it is also extremely disheartening to see the extent to which mathematicians are willing to cozy up to OpenAI and Anthropic, but not to organizers in working class, rural, and/or Black and Brown communities who are fighting to preserve their water and air. 

Before getting too AI-pilled I encourage younger researchers to look into the ways in which AI infrastructure is impacting the lives of your neighbors, friends, and community members. It may also be possible that you are not (yet) directly connected to anyone who is experiencing the brunt of adverse health problems stemming from the sprawling AI industrial complex and in that case I recommend doing some research, perhaps with your cohort, about these issues. A good place to start for the North American context might be \href{https://truthout.org/audio/you-cant-drink-data-native-communities-are-resisting-the-ai-death-cycle/}{this interview with indigenous organizers Krystal Two Bulls and Joseph White Eyes} \cite{Interview}, entitled \textit{You Can’t Drink Data. }

\section{There are many ways to augment your work with this new technology and they are not necessarily created equally from an ethical, moral, or scientific perspective. 
}

To conclude, I’d like to remind you to do your absolute best to distinguish between the abstract potential for AI to be useful for mathematics in a hypothetical reality that might be very different from our own, and the actual potential of the most prominent and popularly used AI models in our world as it currently exists. It is logically consistent for it to be possible for an AI system to be helpful for developing mathematics, and for it to be the case that AI corporations don’t give a shit about mathematics as a discipline and will happily see it burn to the ground if it means a better return on their next quarterly report.

For this reason, I encourage you to think as carefully as possible about which AI systems – if any! – you choose to employ in your research. I hope that if AI is in fact here to stay, the mathematics community will take it upon itself, as it has done in the past, to develop freely accessible and open source models. Perhaps these would be \textit{small} language models with far fewer parameters to tune, such that a researcher can use it to literature search or explore the space of complex structures on the six sphere but that no one can use it to generate a picture of Oprah Winfrey’s head on Elon Musk’s body or whatever other nonsense upon which millions of joules are spent.

I know I don’t have the answers. I am also figuring this out as I go. What I do have is faith in the next generation of mathematicians. And finally, you are not alone. As disappointing as it might be to see senior mathematicians acting in some of the ways that they are, know that there are also older people who genuinely want to see you succeed and who will continue to do what they can to make that happen. You can feel free to reach out to me at the listed email if you’d like to chat more about this.

\bibliographystyle{plainurl}
\bibliography{AIStudents}

\end{document}